\documentclass{amsart}
\usepackage{enumerate, amsmath}

\newtheorem*{thmA}{Theorem A}
\newtheorem{thm}{Theorem}
\newtheorem{lem}[thm]{Lemma}

\theoremstyle{definition}

\theoremstyle{remark}

\numberwithin{equation}{section}

\begin{document}

\title{A quadratic form of $p=3k+1$ primes}

\author{Bat-Od Battseren}
\author{Bayarmagnai Gombodorj}

\email{batoddd@gmail.com}
\thanks{}


\date{\today}


\begin{abstract}
We use Zagier's one-sentence proof approach to show that a prime number $p$ admits a form $p=a^2+ab+b^2$ for some integers, $a$ and $b$, if and only if $p=3$ or $p\equiv 1 \pmod{3}$.

\end{abstract}

\maketitle


\section{Introduction}
In \cite{MR1041893}, Zagier gave a one-sentence proof for the following theorem.
\begin{thm}[Fermat]
Let $p$ be a prime number. There exist integers $x$ and $y$ such that $p=x^2+4y^2$ if and only if $p\equiv 1 \pmod{4}$.
\end{thm}
The main tool was the following lemma.
\begin{lem}
Let $q$ be a prime number, and let $G$ be a finite $q$-group acting on a finite set $S$. Then $|S|\equiv |S^G| \pmod{q}$, where $S^G$ is the set of the  elements of $S$ fixed by the $G$-action.
\end{lem} Zagier considered two different involutive transformations $g_1$ and $g_2$ on the set 
\begin{align*}
S=\{(x,y,z)\in\mathbb{N}^3 : x^2+4yz=p\}.
\end{align*} The transformation $g_1$ sends $(x,y,z)\in S$ to $g_1(x,y,z) = (x,z,y)$. The transformation $g_2$ is defined as a piece-wise linear transformation on $S$. These transformations generate two groups $G_1$ and $G_2$ of two elements acting on $S$. It turns out that the set $S^{G_2}$  has odd number of elements, thus, by Lemma 2,  the set $S^{G_1}$ is non-empty, proving Theorem 1. We use a similar idea to prove the following.
\begin{thmA}
Let $p$ be a prime number. There exist integers $a$ and $b$ such that $p=a^2+ab+b^2$ if and only if $p=3$ or $p\equiv 1 \pmod{3}$.
\end{thmA}
The following theorem will be used.
\begin{thm}[\cite{MR2744747}]
A prime number $p$ satisfies $p\equiv 1 \pmod{3}$ if and only if there exist unique positive integers $x, y \in \mathbb{N}$ such that $p=x^2+3y^2$.
\end{thm}

\section{Proof of Theorem A}
In this section, we prove Theorem A.
The only if part is easily checked by exhausting all possibilities of $a$ and $b$ modulo 3. Let us prove the converse. The case $p=3$ is trivial as $p=1^2+1\cdot 1 + 1^2$. Now, assume $p>3$ and $p\equiv 1 \pmod 3$.
Define 
\begin{align*}
S = \{(x,y,z)\in\mathbb{N}^3:3x^2 +yz=p\}.
\end{align*} Define a transformation $g_0$ on $S$ as $g_0 (x,y,z) = (x,z,y)$. Note that $G_0 = \langle g_0 \rangle_2$ has two elements and acts on $S$. Then, Theorem 3 says that the set $S^{G_0}$ has exactly one element. 

We want to construct another involutive action on $S$. Partition $S$ into the subsets
\begin{align*}
&S_1 = \{ (x,y,z) \in S:6x\leq y-3z\},\\
&S_2= \{ (x,y,z) \in S:5(y-3z)<30x\leq 6(y-2z)\},\\
&S_3 = \{ (x,y,z) \in S:12(y-2z)<60x\leq 5(3x-4z)\},\\
&S_4 = \{ (x,y,z) \in S:7(3y-4z)<84x\leq 24(y-z)\},\\
&S_5 = \{ (x,y,z) \in S:24(y-z)<84x\leq 7(4y-3z)\},\\
&S_6 = \{ (x,y,z) \in S:5(4y-3z)<60x\leq 12(2y-z)\},\\
&S_7 = \{ (x,y,z) \in S:6(2y-z)<30x\leq 5(3y-z)\},\\
&S_8 = \{ (x,y,z) \in S:3y-z<6x\leq 6y\},\\
&S_9 = \{ (x,y,z) \in S:y<x\},\\
&S_{10}= \{ (x,y,z) \in S:y=x\}.
\end{align*}
Define 
\begin{align*}
g (x,y,z) = \left\lbrace \begin{array}{lr}
(x+z,z,-6x+y-3z), & S_1\\
(-5x+y-2z,-12x+3y-4z,6x-y+3z), & S_2\\
(5x-y+2z,6x-y+3z,-12x+3y-4z), & S_3\\
(-7x+2-2z,-12x+4y-3z,12x-3y+4z), & S_4\\
(7x-2y+2z,12x-3y+4z,-12x+4y-3z), & S_5\\
(-5x+2y-z,-6x+3y-z,12x-4y+3z), & S_6\\
(5x-2y+z,12x-4y+3z,-6x+3y-z), & S_7\\
(-x+y,y,6x-3y-z), & S_8\\
(x-y,6x-3y+z,y), & S_9\\
(x,y,z), & S_{10}
\end{array}\right.
\end{align*} Then $g$ is a well defined transformation on $S$. For example, if $(x,y,z)\in S_1$, then we have $3x^2+yz=p$ and $6x\leq y-3z$. Then for $(a,b,c) = g(x,y,z) = (x+z,z,-6x+y-3z)$, we have
$3a^2+bc = 3(x+z)^2 + z (-6x+y-3z) = 3x^2+yz = p$ and $(x+z,z,-6x+y-3z)\in \mathbb{N}^3$. The other cases can be checked in a similar way.

If we denote 
\begin{align*}
B=\begin{pmatrix}
-1 & 1 & 0\\
0 &1 &0\\
6&-3&1
\end{pmatrix},\quad 
X= \begin{pmatrix}
-1 & 0 & 0\\
0 &0 &1\\
0&1&0
\end{pmatrix},\quad 
A=BX=\begin{pmatrix}
1 & 0 & 1\\
0 &0 &1\\
-6&1&-3
\end{pmatrix}
\end{align*}
the transformation $g$ coincides with $T_1 = A$, $T_2=-XA^2$, $T_3 = -A^2$, $T_4 = XA^3$, $T_5 = A^3$, $T_6=-XA^4$, $T_7 = -A^4$, $T_8 = XA^5$, $T_9 = A^5$, and $T_{10} = A^6=I_3$ when restricted on $S_1$, $\dots$, $S_9$, and $S_{10}$, respectively. In other words, for all $1\leq i \leq 10$, we have $g|_{S_i} = T_i$. It is easily checked that 
\begin{align*}\label{formula}
T_1 T_9=T_3 T_7=T_2^2=T_4^2=T_5^2=T_6^2=T_8^2 = T_{10}=I_3.
\end{align*} These identities prove that $g$ is an involution. For example, if $(x,y,z)\in S_9$ then $g(x,y,z)=(x-y,6x-3y+z,y)$ is in $S_1$ since $6(x-y)\leq (6x-3y+z) - 3y$, and $g(g(x,y,z)) = T_1T_9 (x,y,z) = (x,y,z)$. The other cases are handled similarly. Thus, $G = \langle g\rangle_2$ acts on $S$.

From Lemma 2 and Theorem 3, we know the set $S^G$ has odd number of elements. Let us look more into the equation 
\begin{align}
g(x,y,z) = (x,y,z) \text{ on } S.
\end{align} Suppose $(x,y,z)\in S$ is a solution to $(2.1)$. 
\begin{enumerate}
\item If  $(x,y,z)\in S_1$, then $x=x+z>x$, which is impossible.
\item If $(x,y,z)\in S_2$, then $y=6x+2z$ and 
\begin{align*}
p=3x^2 +(6x+2z)z\equiv 2z^2 \not\equiv 1 \pmod{3},
\end{align*} which is impossible.
\item If  $(x,y,z)\in S_3$, then $z=(-12x+3y)/5$ and 
\begin{align*}
p = 3x^2 +y(-12x+3y)/5 \equiv 0 \pmod 3,
\end{align*} which is impossible.
\item If $(x,y,z)\in S_4$, then $y=4x+z$ and 
\begin{align*}
p=3x^2 +(4x+z)z = (x+z)(3x+z),
\end{align*} which is impossible.
\item We skip the case $(x,y,z)\in S_5$ intentionally.
\item If $(x,y,z)\in S_6$, then $y=3x+z/2$ and 
\begin{align*}
p=3x^2 +(3x+z_0)2z_0\equiv 2z_0^2 \not\equiv 1 \pmod{3},
\end{align*} which is impossible.
\item If $(x,y,z)\in S_8$, then $y=2x$ and 
\begin{align*}
p=3x^2 +2xz = x(3x+2z),
\end{align*} so $(x,y,z) = (1,2,\frac{p-3}{2})$.
\item If $(x,y,z)\in S_{10}$, then $y=x$ and 
\begin{align*}
p=3x^2 +xz = x(3x+z),
\end{align*} so $(x,y,z) = (1,1,{p-3})$.
\end{enumerate}
So far, we have exactly 2 solutions. Thus there are odd number of solutions in the set  $S_5$. In particular, there exists an element $(x,y,z)\in S_5 \cap S^G$ which satisfies
\begin{align*}p=3x^2 +yz \quad \text{and}\quad
(x,y,z) = (7x-2y+2z,12x-3y+4z,-12x+4y-3z).
\end{align*}  It follows that $y=3x+z$ and  
\begin{align}
p = 3x^2 +yz = 3x^2 +(3x+z)z = 3x^2 + 3xz + z^2 = x^2 + x(x+z)+ (x+z)^2
\end{align} Now, it is enough to choose $a=x$ and $b=x+z$ to have $p = a^2+ab+b^2$.

In the following table, we provide the first few  examples.
\begin{table}[!htbp]
\begin{tabular}{lllllllllllllllll}
\hline
\multicolumn{1}{|l|}{$p$} & \multicolumn{1}{l|}{7} & \multicolumn{1}{l|}{13} & \multicolumn{1}{l|}{19} & \multicolumn{1}{l|}{31} & \multicolumn{1}{l|}{37} & \multicolumn{1}{l|}{43} & \multicolumn{1}{l|}{61} & \multicolumn{1}{l|}{67} & \multicolumn{1}{l|}{73} & \multicolumn{1}{l|}{79} & \multicolumn{1}{l|}{97} & \multicolumn{1}{l|}{103} & \multicolumn{1}{l|}{109} & \multicolumn{1}{l|}{127} & \multicolumn{1}{l|}{139} & \multicolumn{1}{l|}{151} \\ \hline
\multicolumn{1}{|l|}{$a$} & \multicolumn{1}{l|}{1} & \multicolumn{1}{l|}{1}  & \multicolumn{1}{l|}{2}  & \multicolumn{1}{l|}{1}  & \multicolumn{1}{l|}{3}  & \multicolumn{1}{l|}{1}  & \multicolumn{1}{l|}{4}  & \multicolumn{1}{l|}{2}  & \multicolumn{1}{l|}{1}  & \multicolumn{1}{l|}{3}  & \multicolumn{1}{l|}{3}  & \multicolumn{1}{l|}{2}   & \multicolumn{1}{l|}{5}   & \multicolumn{1}{l|}{6}   & \multicolumn{1}{l|}{3}   & \multicolumn{1}{l|}{5}   \\ \hline
\multicolumn{1}{|l|}{$b$} & \multicolumn{1}{l|}{2} & \multicolumn{1}{l|}{3}  & \multicolumn{1}{l|}{3}  & \multicolumn{1}{l|}{5}  & \multicolumn{1}{l|}{4}  & \multicolumn{1}{l|}{6}  & \multicolumn{1}{l|}{5}  & \multicolumn{1}{l|}{7}  & \multicolumn{1}{l|}{8}  & \multicolumn{1}{l|}{7}  & \multicolumn{1}{l|}{8}  & \multicolumn{1}{l|}{9}   & \multicolumn{1}{l|}{7}   & \multicolumn{1}{l|}{7}   & \multicolumn{1}{l|}{10}  & \multicolumn{1}{l|}{9}   \\ \hline
                        &                        &                         &                         &                         &                         &                         &                         &                         &                         &                         &                         &                          &                          &                          &                          &                          \\
                        &                        &                         &                         &                         &                         &                         &                         &                         &                         &                         &                         &                          &                          &                          &                          &                         
\end{tabular}
\end{table}

\bibliographystyle{amsplain}
\bibliographystyle{alpha}
\bibliography{mybibfile}

\end{document}